# STOCHASTIC INTEGRAL REPRESENTATION AND REGULARITY OF THE DENSITY FOR THE EXIT MEASURE OF SUPER-BROWNIAN MOTION

By Jean-François Le Gall and Leonid Mytnik[1]

*Ecole Normale Supérieure and Technion—Israel Institute of Technology*

This paper studies the regularity properties of the density of the exit measure for super-Brownian motion with $(1+\beta)$-stable branching mechanism. It establishes the continuity of the density in dimension $d=2$ and the unboundedness of the density in all other dimensions where the density exists. An alternative description of the exit measure and its density is also given via a stochastic integral representation. Results are applied to the probabilistic representation of nonnegative solutions of the partial differential equation $\Delta u = u^{1+\beta}$.

**1. Introduction and statement of results.** This paper is devoted to regularity results for the density of the exit measure of super-Brownian motion with $(1+\beta)$-stable branching mechanism from a smooth domain of $\mathbb{R}^d$. Exit measures of superprocesses were introduced by Dynkin in connection with applications to partial differential equations (see in particular [5] and [6]). Here we use a stochastic integral representation of exit measures to get precise information on their regularity or irregularity. As an application we provide a probabilistic representation for all nonnegative solutions of $\Delta u = u^{1+\beta}$ in a smooth domain, in the so-called subcritical case where $d < 1 + 2/\beta$.

Let $D$ be a bounded domain of class $C^2$ in $\mathbb{R}^d$ ($d \geq 2$ throughout this work). If $x \in D$, we write $\rho(x) = \text{dist}(x, D^c)$ for the distance of $x$ to the complement of $D$. We denote by $\mathcal{M}_F^D$ the space of all finite measures on $D$, which is equipped with the weak topology. If $\mu \in \mathcal{M}_F^D$, $\text{supp}(\mu)$ denotes the closed support of $\mu$, which is a subset of $\overline{D}$, and we set

$$\mathcal{M}_{F,c}^D \equiv \{\mu \in \mathcal{M}_F^D : \text{supp}(\mu) \subset D\}.$$

Received June 2003; revised October 2003.

[1]Supported in part by U.S.–Israel Binational Science Foundation Grant 2000065 and Israel Science Foundation Grant 116/01-10.0.

The integral of a function $\phi$ with respect to a measure $\mu$ will often be written as $\langle \mu, \phi \rangle$.

Let $\beta \in (0,1]$ and let $X = (X_t, t \geq 0)$ be a super-Brownian motion in $D$ with $(1+\beta)$-stable branching mechanism. To be specific, $X$ is a superprocess with branching mechanism $\psi(u) = u^{1+\beta}$, whose underlying spatial motion is Brownian motion in $\mathbb{R}^d$ killed when it exits $D$. The process $X$ is a strong Markov process with values in $\mathcal{M}_F^D$, whose distribution will be characterized in Section 2. If $\mu \in \mathcal{M}_F^D$, we write $\mathbb{P}_\mu$ for the probability measure under which $X$ starts from $\mu$. In the first two theorems below, we will consider the case where the initial value $\mu \in \mathcal{M}_{F,c}^D$ (see, however, Remark 3.1).

As a special case of the martingale problem recalled in Section 2.3, we know that for every twice continuously differentiable function $\phi$ on $D$, with compact support contained in $D$,

$$\langle X_t, \phi \rangle = \langle \mu, \phi \rangle + \int_0^t \langle X_s, \tfrac{1}{2}\phi \rangle \, ds + M_t(\phi),$$

where $M_t(\phi)$ is a martingale under $\mathbb{P}_\mu$. It will be convenient to use the notation

$$M_t(\phi) = \int_0^\infty \int_D \mathbb{1}_{[0,t]}(s)\phi(x) M(ds, dx).$$

Standard arguments then show that the "stochastic integral"

$$\int_0^\infty \int_D f(s,x) \, M(ds, dx)$$

can be defined for a wide class of integrands $f$ (see Section 2.3 and the beginning of Section 3).

Let $X^D$ be the exit measure of $X$ from $D$. Note that the usual definition of $X^D$ involves the associated historical process, which contains more information than $(X_t, t \geq 0)$. Alternatively, one can proceed as in [6] or [8] by defining the superprocess as the collection of all exit measures from time-space open sets (these include the measures $X_t$ as special cases). The measure $X^D$ is a random finite measure supported on $\partial D$. We prove in Section 2 that $X^D$ can be obtained via the following approximation, which is of independent interest. For every $\varepsilon > 0$, set

$$D_\varepsilon = \{x \in D : \rho(x) > \varepsilon\}, \qquad F_\varepsilon = D \setminus D_\varepsilon,$$

and

$$X_\varepsilon^D(dy) = \varepsilon^{-2} \int_0^\infty \mathbb{1}_{F_\varepsilon}(y) X_t(dy) \, dt.$$

Then $X_\varepsilon^D$ converges weakly to $X^D$ as $\varepsilon$ tends to 0, in $\mathbb{P}_\mu$-probability (see Proposition 2.1). This shows in particular that $X^D$ is a measurable function of $(X_t, t \geq 0)$.

It was proved in [1] and [21] that $X^D$ is almost surely absolutely continuous with respect to Lebesgue measure on $\partial D$ if and only if $2 \leq d < 1 + 2/\beta$. In the case $\beta = 1$ and $d = 2$, more can be said: $X^D$ has a continuous density (see [15]). In this work we consider the "stable branching" case; that is, from now on we concentrate on the case $0 < \beta < 1$, and we address the question of regularity of the density of the exit measure in dimensions $2 \leq d < 1 + 2/\beta$.

Our first theorem provides a stochastic integral representation for the exit measure and for its density when it exists.

THEOREM 1.1. *Let $(P_D(x,y), x \in D, y \in \partial D)$ denote the Poisson kernel of $D$, and let $\sigma$ denote Lebesgue measure on $\partial D$. Let $\mu \in \mathcal{M}_{F,c}^D$.*

(i) *For every continuous function $\phi$ on $\partial D$, $\mathbb{P}_\mu$-a.s.,*

$$\langle X^D, \phi \rangle = \langle \mu, P_D \phi \rangle + \int_0^\infty \int_D P_D \phi(x) M(ds, dx), \tag{1.1}$$

*where*

$$P_D \phi(x) = \int_{\partial D} P_D(x,y) \phi(y) \sigma(dy).$$

(ii) *Suppose that $d < 1 + 2/\beta$. Then, for every $y \in \partial D$, we may define under $\mathbb{P}_\mu$,*

$$\overline{X}^D(y) = \int_D P_D(x,y) \mu(dx) + \int_0^\infty \int_D P_D(x,y) M(ds, dx). \tag{1.2}$$

*The mapping $y \to \overline{X}^D(y)$ is continuous in $L^p(\mathbb{P}_\mu)$, for any $p \in [1, 1+\beta)$, and we have $\overline{X}^D(y) \geq 0$, $\mathbb{P}_\mu$-a.s., for every $y \in \partial D$. Finally,*

$$X^D(dy) = \overline{X}^D(y) \sigma(dy), \qquad \mathbb{P}_\mu\text{-a.s.}$$

To be precise, we should say in the last assertion that we consider a measurable modification of the process $(\overline{X}^D(y), y \in \partial D)$.

We now come to the main result of the present work, which deals with the regularity properties of $\overline{X}^D$. For any measurable function $f: \partial D \mapsto \mathbb{R}$, let $\|f\|_B$ denote the essential supremum (with respect to Lebesgue measure on $\partial D$) of $f$ on the relative open set $B \subset \partial D$.

THEOREM 1.2 (Regularity and irregularity of density). *Let $\mu \in \mathcal{M}_{F,c}^D$.*

(a) *If $d = 2$, the process $(\overline{X}^D(x), x \in \partial D)$ has a continuous modification under $\mathbb{P}_\mu$.*

(b) *Suppose that $3 \leq d < 1 + 2/\beta$. Then*

$$\|\overline{X}^D(\cdot)\|_U = \infty \qquad \textit{whenever } X^D(U) > 0, \textit{ for any open set } U \subset \partial D, \mathbb{P}_\mu\textit{-a.s.}$$

Obviously, the second part of the theorem remains valid if we replace $\overline{X}^D$ by any version of the Radon–Nikodym derivative of $X^D$ with respect to $\sigma$. Thus, when $3 \leq d < 1 + 2/\beta$, there exists no continuous density of the exit measure.

The main motivation for studying exit measures comes from their connections with partial differential equations. A basic result of Dynkin [5] shows that the exit measure yields a probabilistic solution of the nonlinear Dirichlet problem associated with $\frac{1}{2}\Delta u = u^{1+\beta}$. To be specific, for any nonnegative continuous function $\phi$ on $\partial D$, the function

$$(1.3) \qquad v(x) = -\log \mathbb{E}_{\delta_x}[e^{-\langle X^D, \phi \rangle}], \qquad x \in D,$$

is the unique nonnegative solution to the following boundary value problem in $D$:

$$(1.4) \qquad \begin{aligned} \tfrac{1}{2}\Delta v &= v^{1+\beta} \qquad \text{in } D, \\ v &= \phi \qquad \text{on } \partial D. \end{aligned}$$

A major problem is to extend this probabilistic representation to all nonnegative solutions of $\frac{1}{2}\Delta u = u^{1+\beta}$ in $D$, and to see that this representation induces a one-to-one correspondence between solutions and their traces on the boundary (defined in a proper way). This problem was solved in [15] in the particular case $\beta = 1$, $d = 2$. Later, Marcus and Véron [17] generalized the results of [15] by showing that in the so-called subcritical case $d < 1 + 2/\beta$, there is a one-to-one correspondence between nonnegative solutions and admissible traces. The next theorem gives a probabilistic formula for this correspondence. In order to be able to use the results of [17], we restrict our attention to the case of the unit ball.

We need one more definition. The range $\mathcal{R}$ of $X$ is defined as the closure of the set

$$\bigcup_{t \geq 0} \operatorname{supp}(X_t).$$

THEOREM 1.3. *Suppose that $d < 1 + 2/\beta$ and that $D$ is the unit ball of $\mathbb{R}^d$. Let $K$ be a compact subset of $\partial D$, and let $\nu$ be a Radon measure on $\partial D \setminus K$. The function*

$$(1.5) \ u(x) = -\log \mathbb{E}_{\delta_x}\bigg[\mathbb{1}_{\{\mathcal{R} \cap K = \varnothing\}} \exp\bigg(-\int \overline{X}^D(y)\nu(dy)\bigg)\bigg], \qquad x \in D,$$

*solves the equation $\frac{1}{2}\Delta u = u^{1+\beta}$ in $D$. Conversely, if $u$ is any nonnegative solution of $\frac{1}{2}\Delta u = u^{1+\beta}$ in $D$, there exists a unique pair $(K, \nu)$ such that the representation formula* (1.5) *holds.*

As the proof will show, the pair $(K,\nu)$ can be interpreted as the trace of the solution $u$ (defined analytically in [17]).

Let us emphasize an important point. To make sense of the probabilistic representation stated in Theorem 1.3, it is crucial to have chosen a specified version of the Radon–Nikodym density of the exit measure. In dimension $d=2$, we may of course choose the continuous density (as was done in [15]), but Theorem 1.3 shows that in higher dimensions the right choice is to consider the process $\overline{X}^D(y)$ as defined in Theorem 1.1.

REMARK 1.1. In the present work, we do not discuss the quadratic branching case $\beta=1$. However, our results also hold in that case. Both Theorem 1.2(a) and Theorem 1.3 are proved in [15] in the case $\beta=1$. Furthermore, the reader will easily check that the stochastic integral representation of Theorem 1.1 is also valid in that case: $M$ should then be interpreted as the usual $L^2$-martingale measure associated with super-Brownian motion. As a matter of fact, this stochastic representation can be used to simplify the proof of the key technical lemma of [15].

Let us record some convenient notation for future use. In general, if $F$ is a set of functions, we write $F^+$ for the set of all nonnegative functions in $F$. We use $c$ or $C$ to denote a positive, finite constant whose value may vary from place to place. A notation of the form $c(a,b,\ldots)$ means that this constant depends on parameters $a,b,\ldots$. If $E$ is a metric space, let $\mathcal{B}(E)$ be the corresponding Borel $\sigma$-algebra [$\mathcal{B}(E)$ will also serve as the set of Borel measurable functions on $E$]. We denote by $\mathcal{C}(E)$ the space of all continuous functions on $E$ and by $\mathcal{C}_{\mathrm{b}}(E)$ [resp. $\mathcal{B}_{\mathrm{b}}(E)$] the space of bounded functions in $\mathcal{C}(E)$ [resp. in $\mathcal{B}(E)$]. We also denote by $\mathcal{C}_0^2(D)$ the set of all twice continuously differentiable functions on $D$ with compact support contained in $D$. Finally, if $x \in \mathbb{R}^d$ and $r>0$, $B(x,r)$ stands for the open ball of radius $r$ centered at $x$.

The paper is organized as follows. Section 2 recalls basic facts about super-Brownian motion and states some preliminary results. Theorem 1.1 is proved in Section 3, Theorem 1.2(a) is proved in Section 4, and part (b) of Theorem 1.2 is proved in Section 5. Connections with partial differential equations are discussed in Section 6. The Appendix gives the proof of a technical auxiliary lemma.

## 2. Preliminaries.

2.1. *Estimates for the Green function and the Poisson kernel.* Let $(G_D(x,y); x,y \in D)$ be the Green function of $D$ and recall that $(P_D(x,z); x \in D, z \in \partial D)$ denotes its Poisson kernel. The functions $G_D$ and $P_D$ are continuous on $D \times D$ and $D \times \partial D$, respectively, and they have the following probabilistic interpretation. Let $(\xi_t, t \geq 0; \Pi_x, x \in D)$ denote Brownian motion killed at its first exit time from $D$, and let $\zeta$ be the lifetime of this process. Then, for any $\phi \in \mathcal{B}_b(\mathbb{R}^d)$ and $x \in D$,

$$\Pi_x\left[\int_0^\zeta \phi(\xi_t)\,dt\right] = \int_D \phi(y)G_D(x,y)\,dy,$$

$$\Pi_x[\phi(\xi_{\zeta-})] = \int_{\partial D} \phi(y)P_D(x,y)\sigma(dy).$$

We will use the following estimates. For every $x, y \in D$ and $z \in \partial D$,

(2.1) $$G_D(x,y) \leq C(D)\rho(y)|x-y|^{1-d},$$

(2.2) $$G_D(x,y) \leq C(D)\rho(x)\rho(y)|x-y|^{-d}$$

and

(2.3) $$P_D(x,z) \leq C(D)\rho(x)|x-z|^{-d}.$$

Estimates (2.1) and (2.2) can be found in Theorem 2.3 of [22] in dimension $d \geq 3$. In dimension $d = 2$, they both follow from the more precise bound in Theorem 6.23 of [2]. Finally, (2.3) is a consequence of (2.2) and the interpretation of the Poisson kernel as half the normal derivative of the Green function at the boundary (see Proposition 5.13 in [2]).

2.2. *Super-Brownian motion and its exit measure.* In this section we recall the basic facts about super-Brownian motion that will be used in the proofs of our results, and we also discuss properties of the associated martingale measure. Without additional effort, the results of this section are valid in a more general setting than in the Introduction, namely, for a branching mechanism function $\psi$ of the type

$$\psi(u) = \int n(dr)(e^{-ur} - 1 + ur), \qquad u \geq 0,$$

where $n(dr)$ is a $\sigma$-finite measure on $(0, \infty)$ such that $\int (r \wedge r^2)n(dr) < \infty$. Note that $\psi(u) \geq 0$ for every $u \geq 0$.

Our super-Brownian motion $X$ with branching mechanism $\psi$ is a time-homogeneous Markov process in $\mathcal{M}_F^D$, whose semigroup is characterized as follows: For every $\mu \in \mathcal{M}_F^D, \phi \in \mathcal{B}_b^+(D)$ and $t \geq 0$,

$$\mathbb{E}_\mu[e^{-\langle X_t, \phi \rangle}] = \exp -\langle \mu, u_t \rangle,$$

where the function $(u_t(x), t \geq 0, x \in D)$ is the unique nonnegative solution of the integral equation

$$u_t(x) + \Pi_x\left[\int_0^{t \wedge \zeta} \psi(u_{t-s}(\xi_s))\,ds\right] = \Pi_x[\phi(\xi_t)\mathbb{1}_{\{t < \zeta\}}]$$

(see, e.g., [6] or Chapter 2 of [16]). In particular, taking $\phi = \lambda > 0$, we get an expression for the Laplace transform of $\langle X_t, 1 \rangle$, from which one easily sees that $\mathbb{E}_\mu[\langle X_t, 1 \rangle] \leq \langle \mu, 1 \rangle$ for every $t \geq 0$.

From the preceding Laplace functional, it is not hard to derive that for any $\mu \in \mathcal{M}_F^D, \phi \in \mathcal{B}_b^+(D)$,

$$\mathbb{E}_\mu\left[\exp\left(-\int_0^\infty \langle X_t, \phi \rangle \, dt\right)\right] = \exp - \langle \mu, v \rangle,$$

where the function $(v(x), x \in D)$ is nonnegative and solves the integral equation

$$v(x) + \Pi_x\left[\int_0^\zeta \psi(v(\xi_t)) \, dt\right] = \Pi_x\left[\int_0^\zeta \phi(\xi_t) \, dt\right].$$

In view of approximating the exit measure $X^D$, we now write the following joint Laplace transform. For any $g \in \mathcal{B}_b^+(\partial D)$ and $\phi \in \mathcal{B}_b^+(D)$,

$$\mathbb{E}_\mu\left[\exp\left(-\int_0^\infty \langle X_t, \phi \rangle \, dt - \langle X^D, g \rangle\right)\right] = \exp - \langle \mu, w \rangle,$$

where the function $(w(x), x \in D)$ is nonnegative and solves the integral equation

$$w(x) + \Pi_x\left[\int_0^\zeta \psi(w(\xi_t)) \, dt\right] = \Pi_x\left[\int_0^\zeta \phi(\xi_t) \, dt + g(\xi_{\zeta-})\right].$$

This statement is a special case of Theorem I.1.8 in [6]. We can now prove the approximation of the exit measure stated in the Introduction.

PROPOSITION 2.1. *Let $X_\varepsilon^D$ be defined as in Section 1. Then $X_\varepsilon^D$ converges weakly to $X^D$ as $\varepsilon \downarrow 0$, in $\mathbb{P}_\mu$-probability.*

PROOF. Let $\varphi \in \mathcal{C}^+(\overline{D})$. It is enough to prove that

$$\langle X_\varepsilon^D, \varphi \rangle \longrightarrow \langle X^D, \varphi \rangle$$

in $\mathbb{P}_\mu$-probability, as $\varepsilon \to 0$. To this end, we need only check that, for every $\lambda, \lambda' \geq 0$,

$$\mathbb{E}_\mu[\exp(-\lambda \langle X_\varepsilon^D, \varphi \rangle - \lambda' \langle X^D, \varphi \rangle)] \to \mathbb{E}_\mu[\exp(-(\lambda + \lambda')\langle X^D, \varphi \rangle)]$$

as $\varepsilon \to 0$. We fix $\lambda$ and $\lambda'$ and establish the preceding limit.

By our definition of $X_\varepsilon^D$, and results recalled before the statement of the proposition, we have

$$\mathbb{E}_\mu[\exp(-\lambda \langle X_\varepsilon^D, \varphi \rangle - \lambda' \langle X^D, \varphi \rangle)] = \exp(-\langle \mu, w^\varepsilon \rangle),$$

where

$$w^\varepsilon(x) + \Pi_x\left[\int_0^\zeta \psi(w^\varepsilon(\xi_t))\,dt\right]$$
(2.4)
$$= \Pi_x\left[\frac{\lambda}{\varepsilon^2}\int_0^\zeta \mathbb{1}_{F_\varepsilon}(\xi_t)\varphi(\xi_t)\,dt + \lambda'\varphi(\xi_{\zeta-})\right] \equiv h^\varepsilon(x).$$

Similarly,

$$\mathbb{E}_\mu[\exp(-(\lambda+\lambda')\langle X^D,\varphi\rangle)] = \exp(-\langle\mu,w\rangle),$$

where

(2.5) $$w(x) + \Pi_x\left[\int_0^\zeta \psi(w(\xi_t))\,dt\right] = \Pi_x[(\lambda+\lambda')\varphi(\xi_{\zeta-})] \equiv h(x).$$

By standard arguments (see, e.g., the proof of Theorem 1.1 in [5]), (2.5) is equivalent to the boundary value problem

$$\tfrac{1}{2}\Delta w = \psi(w) \quad \text{in } D,$$
$$w = (\lambda+\lambda')\varphi \quad \text{on } \partial D.$$

Uniqueness of the nonnegative solution for this boundary value problem is a consequence of the maximum principle, and so we see that $w$ is the unique nonnegative solution of (2.5).

An application of the bounds (2.1) and (2.2) shows that there exists a constant $C(D)$ such that for every $x \in D$ and $\varepsilon \in (0,1]$,

$$\Pi_x\left[\int_0^\zeta \mathbb{1}_{F_\varepsilon}(\xi_t)\,dt\right] = \int_{F_\varepsilon} G_D(x,y)\,dy \leq C(D)\varepsilon^2.$$

To get this, first note that by the strong Markov property it is enough to consider the case when $x \in F_\varepsilon$, and then use the bound (2.1) when $|y-x| \leq \varepsilon$ and the bound (2.2) when $|y-x| > \varepsilon$. The point is to observe that the Lebesgue measure of $F_\varepsilon \cap B(x,\delta)$ is bounded above by $C'(D)\varepsilon\delta^{d-1}$ for every $\delta \in [\varepsilon,\infty)$.

It follows from the previous bound that the functions $h^\varepsilon$, $\varepsilon \in (0,1]$, are uniformly bounded over $D$, and by (2.4) the same holds for the functions $w^\varepsilon$, $\varepsilon \in (0,1]$. We have then

$$\frac{1}{\varepsilon^2}\Pi_x\left[\int_0^\zeta \mathbb{1}_{F_\varepsilon}(\xi_t)\varphi(\xi_t)\,dt\right] = \frac{1}{\varepsilon^2}\int_{F_\varepsilon} G_D(x,y)\varphi(y)\,dy.$$

Using either of the bounds (2.1) or (2.2), and the fact that $P_D(x,z)$ is half the normal derivative of the mapping $y \to G_D(x,y)$ at $z$ [in other words, $G_D(x,y) \sim 2\rho(y)P_D(x,z)$ when $y$ tends to $z$ along the normal to $\partial D$ at $z$], we easily get

$$\lim_{\varepsilon \to 0}\frac{1}{\varepsilon^2}\int_{F_\varepsilon} G_D(x,y)\varphi(y)\,dy = \int_{\partial D} P_D(x,z)\varphi(z)\sigma(dz) = \Pi_x[\varphi(\xi_{\zeta-})].$$

It follows that $h^\varepsilon(x) \to h(x)$ as $\varepsilon \to 0$, for every $x \in D$.

Let $K$ be a compact subset of $D$, and let $\varepsilon_0 \in (0, 1]$ such that $\rho(x) > \varepsilon_0$ for every $x \in K$. Denote by $\zeta_0$ the first exit time from $D_{\varepsilon_0}$. From (2.4) and the strong Markov property at time $\zeta_0$, we get that for every $x \in D_{\varepsilon_0}$ and $\varepsilon \in (0, \varepsilon_0]$,

$$w^\varepsilon(x) + \Pi_x\left[\int_0^{\zeta_0} \psi(w^\varepsilon(\xi_s))\,ds\right] = h_0^\varepsilon(x),$$

where the functions $h_0^\varepsilon$ are harmonic on $D_{\varepsilon_0}$ and uniformly bounded. As previously, this integral equation implies that $w^\varepsilon$ solves $\frac{1}{2}\Delta w^\varepsilon = \psi(w^\varepsilon)$ in $D_{\varepsilon_0}$ and since the functions $w^\varepsilon$ are uniformly bounded on $D$, standard analytic arguments (see, e.g., Theorem 3.9 in [12]) show that the functions $w^\varepsilon$ are equicontinuous on $K$. At least along a subsequence, we may therefore assume that $w^\varepsilon$ converges to a limiting function $\tilde{w}$, uniformly on every compact subset of $D$. By passing to the limit in (2.4), we see that $\tilde{w}$ solves (2.5) and thus $\tilde{w} = w$. We conclude that $w^\varepsilon$ converges to $w$, which completes the proof. $\square$

2.3. *The associated martingale measure.* For the results of this section, it is convenient to equip the underlying probability space $\Omega$ with the filtration $(\mathcal{F}_t)$ generated by $X$, which is completed as usual with the class of $\mathcal{F}_\infty$-measurable sets which are $\mathbb{P}_\mu$-negligible for every $\mu \in \mathcal{M}_F^D$. All martingales or local martingales will be relative to the filtration $(\mathcal{F}_t)$. We will use the standard notation $\Delta X_s = X_s - X_{s-}$ for the jump of $X$ at time $s$ (no confusion should arise from the use of $\Delta$ also for the Laplacian).

We first recall from [3], Section 6.1 or [10] that $X$ satisfies the following martingale problem. For every $\varphi \in \mathcal{C}_0^2(D)$ and every $f \in \mathcal{C}^2(\mathbb{R})$,

$$f(\langle X_t, \varphi\rangle) - f(\langle X_0, \varphi\rangle) - \tfrac{1}{2}\int_0^t f'(\langle X_s, \varphi\rangle)\langle X_s, \Delta\varphi\rangle\,ds$$
$$- \int_0^t \left(\int_D \int_{(0,\infty)} (f(\langle X_s, \varphi\rangle + r\varphi(x))\right.$$
$$\left. - f(\langle X_s, \varphi\rangle) - f'(\langle X_s, \varphi\rangle)r\varphi(x))n(dr)X_s(dx)\right)ds$$

is a local martingale.

From this martingale problem, one easily infers that the jumps of $X$ must be of the following type. If $s > 0$ is a jump time of $X$, then $\Delta X_s = r\delta_x$ for some $r > 0$ and $x \in \mathbb{R}^d$. More precisely, if $J$ denotes the set of all jump times of $X$, the compensator of the random measure

$$N := \sum_{s \in J} \delta_{(s, \Delta X_s)}$$

is given by the following formula. For any nonnegative predictable function $F$ on $\mathbb{R}_+ \times \Omega \times \mathcal{M}_F^D$,

$$(2.6) \qquad \mathbb{E}_\mu\bigg[\sum_{s \in J} F(s, \omega, \Delta X_s)\bigg] = \mathbb{E}_\mu\bigg[\int F(s, \omega, \mu) \widehat{N}(ds, d\mu)\bigg],$$

where $\widehat{N}$ is the random measure on $\mathbb{R}_+ \times \mathcal{M}_F^D$ defined by

$$\int G(s,\mu) \widehat{N}(ds,d\mu) = \int_0^\infty ds \int n(dr) \int X_s(dx) G(s, r\delta_x).$$

See Théorème 7 in [10] or [3], page 111.

Let $F$ be a measurable function on $\mathbb{R}_+ \times \mathcal{M}_F^D$ such that for every $t \geq 0$,

$$(2.7) \qquad \mathbb{E}_\mu\bigg[\bigg(\sum_{s \in J \cap [0,t]} F(s, \Delta X_s)^2\bigg)^{1/2}\bigg] < \infty.$$

Following [14], Section II.1d, we can then define the stochastic integral of $F$ with respect to the compensated measure $N - \widehat{N}$,

$$\int_0^t F(s, \mu)(N - \widehat{N})(ds, d\mu),$$

as the unique purely discontinuous martingale (vanishing at time 0) whose jumps are indistinguishable from the process $\mathbb{1}_J(s) F(s, \Delta X_s)$.

We shall be interested in the special case where $F(s, \mu) = F_\phi(s, \mu) \equiv \int \phi(s, x) \mu(dx)$ for some measurable function $\phi$ on $\mathbb{R}_+ \times D$. [Some convention is needed when $\int |\phi(s,x)| \mu(dx) = \infty$, but this will be irrelevant in what follows.] If $\phi$ is bounded, then it is easy to see that condition (2.7) holds. Indeed, we can bound separately

$$\mathbb{E}_\mu\bigg[\bigg(\sum_{s \leq t} \langle \Delta X_s, 1\rangle^2 \mathbb{1}_{\{\langle \Delta X_s, 1\rangle \leq 1\}}\bigg)^{1/2}\bigg]$$
$$\leq \mathbb{E}_\mu\bigg[\sum_{s \leq t} \langle \Delta X_s, 1\rangle^2 \mathbb{1}_{\{\langle \Delta X_s, 1\rangle \leq 1\}}\bigg]^{1/2}$$
$$= \bigg(\int_{(0,1]} r^2 n(dr) \mathbb{E}_\mu\bigg[\int_0^t \langle X_s, 1\rangle \, ds\bigg]\bigg)^{1/2} < \infty,$$

and, using the simple inequality $a_1^2 + \cdots + a_n^2 \leq (a_1 + \cdots + a_n)^2$ for any nonnegative reals $a_1, \ldots, a_n$,

$$\mathbb{E}_\mu\bigg[\bigg(\sum_{s \leq t} \langle \Delta X_s, 1\rangle^2 \mathbb{1}_{\{\langle \Delta X_s, 1\rangle > 1\}}\bigg)^{1/2}\bigg] \leq \mathbb{E}_\mu\bigg[\sum_{s \leq t} \langle \Delta X_s, 1\rangle \mathbb{1}_{\{\langle \Delta X_s, 1\rangle > 1\}}\bigg]$$
$$= \int_{(1,\infty)} r n(dr) \mathbb{E}_\mu\bigg[\int_0^t \langle X_s, 1\rangle \, ds\bigg] < \infty.$$

In both cases, we have used (2.6) and the fact that $\mathbb{E}_\mu[\langle X_t, 1\rangle] \leq \langle \mu, 1\rangle$.

To simplify notation, we write

$$M_t(\phi) = \int_0^t \int_D \phi(s,x) M(ds,dx) \equiv \int_0^t F_\phi(s,\mu)(N - \widehat{N})(ds, d\mu),$$

whenever (2.7) holds for $F = F_\phi$. This is consistent with the notation of the Introduction. Indeed, if $\phi(s,x) = \varphi(x)$ where $\varphi \in \mathcal{C}_0^2(D)$, then by the very definition, $M_t(\phi)$ is a purely discontinuous martingale with the same jumps as the process $\langle X_t, \varphi \rangle$. Since the same holds for the process

$$\widetilde{M}_t(\varphi) := \langle X_t, \varphi \rangle - \langle X_0, \varphi \rangle - \tfrac{1}{2} \int_0^t \langle X_s, \Delta \varphi \rangle\, ds$$

(see Théorème 7 in [10]), we get that $M_t(\phi) = \widetilde{M}_t(\varphi)$.

**3. The stochastic integral representation.** We return to the special case where $\psi(u) = u^{1+\beta}$ and thus

$$n(dr) = \frac{\beta(\beta+1)}{\Gamma(1-\beta)} r^{-2-\beta}\, dr$$

for some $\beta \in (0,1)$.

In this section and in the next two, we fix the initial measure $\mu$ of our super-Brownian motion, and we assume that $\mu \in \mathcal{M}_{F,c}^D$. To simplify notation, we write $\mathbb{P}$ instead of $\mathbb{P}_\mu$ and $\mathbb{E}$ instead of $\mathbb{E}_\mu$.

We need to introduce some notation. Let $\{p_t^D(x,y),\ t>0, x,y \in D\}$ be the transition density of Brownian motion killed on its exit from $D$, and let $\{S_t^D,\ t \geq 0\}$ be the corresponding semigroup. For any measure $\nu \in \mathcal{M}_F^D$ set

$$S_t^D \nu(y) = \int p_t^D(x,y) \nu(dx), \qquad y \in D,\ t > 0.$$

Recall that $\mathbb{E}[\langle X_t, \phi \rangle] = \int_D \phi(x) S_t^D \mu(x)\, dx$ for every $t \geq 0$ and $\phi \in \mathcal{B}^+(D)$ (this first-moment formula is easy from the Laplace functional of $\langle X_t, \phi \rangle$ recalled in Section 2).

For any $p \geq 1$, we define the Banach space

$$\mathbb{L}^p \equiv L^p(\mathbb{R}_+ \times D, S_s^D \mu(x)\, ds\, dx)$$

of equivalent classes of measurable functions with finite norms

$$\|f\|_p \equiv \left( \int_0^\infty \int_D |f(s,x)|^p S_s^D \mu(x)\, dx\, ds \right)^{1/p}.$$

Note that if $f$ does not depend on the "time" parameter $s$, then

$$\|f\|_p = \left( \int_D |f(x)|^p G_D \mu(x)\, dx \right)^{1/p},$$

where
$$G_D\nu(y) \equiv \int_D G_D(x,y)\nu(dx), \qquad y \in D, \ \nu \in \mathcal{M}_F^D.$$

LEMMA 3.1. *Let $\phi \in \mathbb{L}^p$, for some $p \in (1+\beta, 2)$. Then the martingale*
$$M_t(\phi) = \int_0^t \int_D \phi(s,x) M(ds, dx), \qquad t \geq 0,$$
*is well defined, and bounded in $L^q(\mathbb{P})$ for every $q \in (1, 1+\beta)$. More precisely, for every $q \in (1, 1+\beta)$,*

(3.1) $$\mathbb{E}\left[\sup_{t \geq 0}|M_t(\phi)|^q\right] \leq c(\beta, p, q)(\|\phi\|_p^q + \|\phi\|_q^q).$$

*Moreover, for any sequence of functions $\{\phi_n, n \geq 1\}$ such that $\phi_n \to \phi$ in $\mathbb{L}^p$, as $n \to \infty$, we have*

(3.2) $$\lim_{n \to \infty} \mathbb{E}\left[\sup_{t \geq 0}|M_t(\phi_n) - M_t(\phi)|^q\right] = 0 \qquad \forall q \in (1, 1+\beta).$$

PROOF. To see that the martingale $M_t(\phi)$ is well defined, we need to verify condition (2.7) with $F = F_\phi$. We will in fact prove more by checking that, for every $q \in (1, 1+\beta)$,

(3.3) $$\mathbb{E}\left[\left(\sum_{s \in J} F(s, \Delta X_s)^2\right)^{q/2}\right] < \infty.$$

First note that since $p/2 \leq 1$, we have $(\sum_{i \in I} a_i)^{p/2} \leq \sum_{i \in I} a_i^{p/2}$ whenever $a_i \geq 0$ for every $i \in I$. We use this in the second inequality below:

$$\mathbb{E}\left[\left(\sum_{s \in J} \mathbb{1}_{\{\langle \Delta X_s, 1 \rangle \leq 1\}} F(s, \Delta X_s)^2\right)^{q/2}\right]$$
$$\leq \mathbb{E}\left[\left(\sum_{s \in J} \mathbb{1}_{\{\langle \Delta X_s, 1 \rangle \leq 1\}} F(s, \Delta X_s)^2\right)^{p/2}\right]^{q/p}$$
$$\leq \mathbb{E}\left[\sum_{s \in J} \mathbb{1}_{\{\langle \Delta X_s, 1 \rangle \leq 1\}} |F(s, \Delta X_s)|^p\right]^{q/p}$$
$$= \mathbb{E}\left[\int_0^\infty ds \int n(dr) \int X_s(dx) \mathbb{1}_{\{r \leq 1\}} r^p |\phi(s,x)|^p\right]^{q/p}$$
$$= \left(\left(\int_{(0,1]} r^p n(dr)\right) \int_0^\infty ds \int_D dx\, S_s^D \mu(x) |\phi(s,x)|^p\right)^{q/p}$$
$$= C(\beta, p, q) \|\phi\|_p^q,$$

using (2.6) and the fact that $\int_{(0,1]} r^p n(dr) < \infty$ since $p > 1 + \beta$.

Similarly,

$$\mathbb{E}\left[\left(\sum_{s \in J} \mathbb{1}_{\{\langle \Delta X_s, 1 \rangle > 1\}} F(s, \Delta X_s)^2\right)^{q/2}\right]$$

$$\leq \mathbb{E}\left[\sum_{s \in J} \mathbb{1}_{\{\langle \Delta X_s, 1 \rangle > 1\}} |F(s, \Delta X_s)|^q\right]$$

$$= \mathbb{E}\left[\int_0^\infty ds \int n(dr) \int X_s(dx) \mathbb{1}_{\{r > 1\}} r^q |\phi(s, x)|^q\right]$$

$$= \left(\left(\int_{(1,\infty)} r^q n(dr)\right) \int_0^\infty ds \int_D S_s^D \mu(dx) |\phi(s, x)|^q\right)$$

$$= C(\beta, q) \|\phi\|_q^q,$$

using (2.6) and the fact that $\int_{(1,\infty)} r^q n(dr) < \infty$ since $q < 1 + \beta$.

By combining the last two bounds, we see that (3.3) holds. Furthermore, by the Burkholder–Davis–Gundy inequality for purely discontinuous martingales (see, e.g., Chapter VII of [4]),

$$\mathbb{E}\left[\sup_{t \geq 0} |M_t(\phi)|^q\right] \leq C(q) \mathbb{E}\left[\left(\sum_{s \in J} F(s, \Delta X_s)^2\right)^{q/2}\right],$$

and the bound (3.1) follows from the previous inequalities. The last assertion is immediate from (3.1), observing that $\phi_n \to \phi$ in $\mathbb{L}^p$ implies $\phi_n \to \phi$ in $\mathbb{L}^q$ since the measure $S_s^D \mu(x) \, dx \, ds$ is finite. $\square$

The next lemma is a Fubini-like theorem for our stochastic integrals.

LEMMA 3.2. *Let $(E, \mathcal{E}, \nu)$ be a $\sigma$-finite measure space and let $\phi$ be a measurable function on $\mathbb{R}_+ \times D \times E$. Assume that for some $p \in (1 + \beta, 2)$,*

$$\int_E \int_0^\infty \int_D |\phi(s, x, y)|^p S_s^D \mu(x) \, dx \, ds \, \nu(dy) < \infty,$$

*and for every $y \in E$,*

$$\int_0^\infty \int_D |\phi(s, x, y)|^p S_s^D \mu(x) \, dx \, ds < \infty.$$

*For every $y \in E$ set $\phi_y(t, x) = \phi(t, x, y)$ and*

$$M_t(\phi_y) = \int_0^t \int_D \phi_y(s, x) M(ds, dx).$$

*Then, for every $t \in [0, \infty]$, the process $(M_t(\phi_y), y \in E)$ has a measurable modification, and*

$$(3.4) \quad \int_E M_t(\phi_y)\nu(dy) = \int_0^t \int_D \left( \int_E \phi(s,x,y)\nu(dy) \right) M(ds,dx), \qquad \mathbb{P}\text{-a.s.}$$

PROOF. We only sketch the arguments. First note that our integrability assumptions guarantee that the stochastic integrals $M_t(\phi_y)$ are well defined for every $y \in E$, that the function $y \to \phi(s,x,y)$ is $\nu$-integrable $S_s^D\mu(x)\,dx\,ds$-a.e., and that the stochastic integral in the right-hand side of (3.4) is well defined, independently of the value we give to $\int_E \phi(s,x,y)\nu(dy)$ when $y \to \phi(s,x,y)$ is not $\nu$-integrable. By standard arguments, it suffices to prove the lemma when $\nu$ is a finite measure and $\phi = \mathbb{1}_A$ is an indicator function (note that the integrability assumptions of the lemma are then automatically satisfied). In the particular case where $A = A_1 \times A_2$, with $A_1 \in \mathcal{B}(\mathbb{R}_+ \times D)$ and $A_2 \in \mathcal{E}$, the various assertions of the lemma are immediately verified. The general case follows from a classical monotone class argument. $\square$

PROOF OF THEOREM 1.1(i). Let $\phi \in \mathcal{C}(\partial D)$. We may extend $\phi$ to a continuous function on $\overline{D}$, which we still denote by $\phi$. By standard techniques (see Proposition 2.13 in [11] or Exercise II.5.2 in [19] for the finite variance branching case), it is easy to obtain that for every $t \geq 0$, $\mathbb{P}$-a.s.,

$$(3.5) \qquad \langle X_t, \phi \rangle = \langle \mu, S_t^D \phi \rangle + \int_0^t \int_D S_{t-s}^D \phi(x) M(ds, dx).$$

We then apply Lemma 3.2 to the (bounded) function $(s,x,t) \to \mathbb{1}_{\{s \leq t\}} S_{t-s}^D \phi(x)$, noting that $\int_0^\infty \mathbb{1}_{\{s \leq t\}} S_{t-s}^D \phi(x)\,dt = G_D \phi(x)$. It follows that

$$(3.6) \qquad \int_0^\infty \langle X_t, \phi \rangle\,dt = \langle \mu, G_D \phi \rangle + \int_0^\infty \int_D G_D \phi(x) M(ds, dx).$$

From the definition of $X_\varepsilon^D$, we get for any $\varepsilon > 0$,

$$(3.7) \qquad \langle X_\varepsilon^D, \phi \rangle = \langle \mu, G_D f^\varepsilon \rangle + \int_0^\infty \int_D G_D f^\varepsilon(x) M(ds, dx),$$

where $f^\varepsilon(x) = \varepsilon^{-2} \mathbb{1}_{F_\varepsilon}(x)\phi(x)$. As in the proof of Proposition 2.1, it is easy to verify that, for every $x \in D$,

$$(3.8) \qquad G_D f^\varepsilon(x) \to P_D \phi(x)$$

as $\varepsilon \to 0$, and furthermore, the functions $G_D f^\varepsilon$ are uniformly bounded over $D$. By dominated convergence, we see that $G_D f^\varepsilon$ converges to $P_D \phi$ in $\mathbb{L}^p$ for every $p \in (1+\beta, 2)$. By passing to the limit $\varepsilon \to 0$ (using the last assertion of Lemma 3.1), we get the desired result. $\square$

PROOF OF THEOREM 1.1(ii). Let $p \in (1, \frac{d+1}{d-1})$. From the bounds (2.1) and (2.3), it is straightforward to verify that, for any compact subset $K$ of $D$,

$$(3.9) \qquad \sup_{x \in K, z \in \partial D} \left( \int_D G_D(x,y) P_D(y,z)^p \, dy \right) < \infty.$$

We are assuming $d < 1 + 2/\beta$, or equivalently $1 + \beta < \frac{d+1}{d-1}$. We can thus choose $p \in (1 + \beta, \frac{d+1}{d-1} \wedge 2)$ and the preceding estimate implies that the (time-independent) function $(s,y) \to P_D(y,z)$ is in $\mathbb{L}^p$ for every $z \in \partial D$. In particular, the stochastic integral appearing in the definition of $\overline{X}^D(z)$ is well defined according to Lemma 3.1. Furthermore, using the fact that $\mu \in \mathcal{M}_{F,c}^D$, we can apply Lemma 3.2 to the function $(y,s,z) \to P_D(y,z)$ and the measurable space $(E, \mathcal{E}, \nu) = (\partial D, \mathcal{B}(\partial D), \sigma)$. It readily follows that the process $(\overline{X}^D(z), z \in \partial D)$ has a measurable modification, and that, for any $\phi \in \mathcal{C}(\partial D)$, $\mathbb{P}$-a.s.,

$$\begin{aligned}
\langle X^D, \phi \rangle &= \langle \mu, P_D \phi \rangle + \int_0^\infty \int_D \left( \int_{\partial D} P_D(x,z) \phi(z) \sigma(dz) \right) M(ds, dx) \\
&= \int_{\partial D} \phi(z) \left( \int_D P_D(x,z) \mu(dx) \right) \sigma(dz) \\
&\quad + \int_{\partial D} \phi(z) \left( \int_0^\infty \int_D P_D(x,z) M(ds, dx) \right) \sigma(dz) \\
&= \int_{\partial D} \phi(z) \overline{X}^D(z) \sigma(dz).
\end{aligned}$$

This is enough to conclude that $X^D(dz) = \overline{X}^D(z) \sigma(dz)$, $\mathbb{P}$-a.s.

In particular, we must have $\overline{X}^D(z) \geq 0$, $\sigma(dz)$ a.e., $\mathbb{P}$-a.s. From the estimate (3.9) and the last assertion of Lemma 3.1, it is easy to see that the mapping $z \to \overline{X}^D(z)$ is continuous in $L^q(\mathbb{P})$, for every $q < 1 + \beta$, and it follows that $\overline{X}^D(z) \geq 0$, $\mathbb{P}$-a.s., for every $z \in \partial D$. □

REMARK 3.1. The proof of part (i) of Theorem 1.1 does not depend on the assumption $\mu \in \mathcal{M}_{F,c}^D$, and the result is indeed true for an initial measure $\mu \in \mathcal{M}_F^D$. Things go differently for part (ii): If $\mu \in \mathcal{M}_F^D \setminus \mathcal{M}_{F,c}^D$, the function $(s,x) \to P_D(x,z)$ may no longer be in $\mathbb{L}^p$ for any $p > 1 + \beta$, and the stochastic integral appearing in (1.2) may not be defined. Still from the additivity property of superprocesses, we can recover from the particular case $\mu \in \mathcal{M}_{F,c}^D$ the fact that the exit measure is absolutely continuous with respect to Lebesgue measure on the boundary.

**4. Continuity of the density in two dimensions.** In this section, we assume that $d = 2$ and we prove part (a) of Theorem 1.2. As we want to use the Riemann mapping theorem, we will first assume that $D$ is simply connected.

The first term in the right-hand side of (1.2) is obviously continuous in $y$. So, to prove the existence of a continuous modification of $\overline{X}^D(y)$, it is enough to check the existence of a continuous modification of the stochastic integral

$$Z(y) \equiv \int_0^\infty \int_D P_D(x,y) M(ds,dx).$$

Before we continue, let us introduce the following notation. Let $D_0$ be the unit disc of the plane and denote by $\sigma_0(dy)$ the Lebesgue measure on the unit circle $\partial D_0$. The Poisson kernel in this case can be computed explicitly:

$$(4.1) \qquad P_0(x,y) = \frac{1}{2\pi} \frac{1-|x|^2}{|y-x|^2}, \qquad x \in D_0,\ y \in \partial D_0.$$

The next lemma is crucial for estimating the moments of increments of $Z(\cdot)$.

LEMMA 4.1. (a) *Set $\rho_0(x) \equiv \mathrm{dist}(x, \partial D_0)$. Let $a \geq 0$, $p \in (0, 2+a)$, and*

$$(4.2) \qquad \gamma = \begin{cases} 2+a-p, & \text{if } \dfrac{2+a}{2} < p < 2+a, \\[4pt] \dfrac{2+a}{2} - \varepsilon, & \text{if } p = \dfrac{2+a}{2}, \\[4pt] p, & \text{if } 0 < p < \dfrac{2+a}{2}, \end{cases}$$

*where $\varepsilon \in (0, \frac{2+a}{2})$ is arbitrary. Then there exists a constant $c = c(p, a, \varepsilon)$ such that*

$$(4.3) \qquad \int_{D_0} \rho_0(x)^a |P_0(x,y_1) - P_0(x,y_2)|^p\, dx \leq c|y_1 - y_2|^\gamma \qquad \forall\, y_1, y_2 \in \partial D_0.$$

(b) *For any $B \subset D_0$ such that $\mathrm{dist}(B, \partial D_0) > 0$, there exists $c = c(B)$ such that*

$$(4.4) \qquad \sup_{x \in B} |P_0(x,y_1) - P_0(x,y_2)| \leq c|y_1 - y_2| \qquad \forall\, y_1, y_2 \in \partial D_0.$$

The proof of Lemma 4.1 appears in the Appendix.

Since $D$ is a bounded simply connected domain in $\mathbb{R}^2$, the Riemann mapping theorem allows us to find a conformal mapping $\psi$ from $D_0$ onto $D$. Under our assumption that $D$ is of class $C^2$, $\psi$ extends to a one-to-one continuous mapping from $\overline{D}_0$ onto $\overline{D}$. In fact, we can say more. According to Chapter 3 of [20], $\psi'$ also has a continuous extension to $\overline{D}_0$ and $\psi'$ does not vanish on $\overline{D}_0$. In particular, $|\psi'|$ is bounded below and above on $\overline{D}_0$ by positive constants. It is also easy to check that for every $x, y \in D$ and $z \in \partial D$,

$$(4.5) \qquad P_D(x,z) = |\psi'(\psi^{-1}(z))|^{-1} P_0(\psi^{-1}(x), \psi^{-1}(z)),$$

$$(4.6) \qquad G_D(x,y) = G_0(\psi^{-1}(x), \psi^{-1}(y)).$$

Let $\theta(x) = |\psi'(\psi^{-1}(x))|$ for every $x \in \overline{D}$.

LEMMA 4.2. *Let $p \in (0,3)$, $\varepsilon \in (0, \frac{3}{2})$ and $\mu \in \mathcal{M}_{F,c}^D$. There exists $c = c(p, \varepsilon, D, \mu)$ such that for every $y_1, y_2 \in \partial D$,*

$$\int_D |\theta(y_1) P_D(x, y_1) - \theta(y_2) P_D(x, y_2)|^p G_D \mu(x) \, dx$$

$$\leq \begin{cases} c|y_1 - y_2|^p, & \text{if } 0 < p < \frac{3}{2}, \\ c|y_1 - y_2|^{3/2 - \varepsilon}, & \text{if } p = \frac{3}{2}, \\ c|y_1 - y_2|^{3-p}, & \text{if } \frac{3}{2} < p < 3. \end{cases}$$

PROOF. Let $\tilde{\mu}$ be the image of $\mu$ under $\psi^{-1}$, and set $\mathcal{E}(\tilde{\mu}) = \text{supp}(\tilde{\mu})$, $\rho_* = \text{dist}(\mathcal{E}(\tilde{\mu}), \partial D_0)$, and $\mathcal{E}(\tilde{\mu})^{\rho_*/2} = \{x \in D_0 : \text{dist}(x, \mathcal{E}(\tilde{\mu})) < \rho_*/2\}$. Then, using (4.5) and (4.6),

$$\int_D |\theta(y_1) P_D(x, y_1) - \theta(y_2) P_D(x, y_2)|^p G_D \mu(x) \, dx$$

$$= \int_{D_0} |P_0(x', \psi^{-1}(y_1)) - P_0(x', \psi^{-1}(y_2))|^p$$

$$\times \left( \int_{D_0} G_0(w, x') \tilde{\mu}(dw) \right) |\psi'(x')|^2 \, dx'.$$

By an application of the Fubini theorem,

$$\int_{D_0} \left( \int_{D_0} G_0(w, x') \tilde{\mu}(dw) \right) |\psi'(x')|^2 \, dx' \leq c(\tilde{\mu}, \psi),$$

and on the other hand, the bounds (2.1) easily imply that for every $x' \in D_0 \setminus \mathcal{E}(\tilde{\mu})^{\rho_*/2}$,

$$\int_{D_0} G_0(w, x') \tilde{\mu}(dw) \leq c(\tilde{\mu}, \psi) \rho_0(x').$$

It readily follows that

$$\int_D |\theta(y_1) P_D(x, y_1) - \theta(y_2) P_D(x, y_2)|^p G_D \mu(x) \, dx$$

$$\leq c(\tilde{\mu}, \psi) \left( \sup_{x' \in \mathcal{E}(\tilde{\mu})^{\rho_*/2}} |P_0(x', \psi^{-1}(y_1)) - P_0(y', \psi^{-1}(y_2))|^p \right.$$

$$\left. + \int_{D_0 \setminus \mathcal{E}(\tilde{\mu})^{\rho_*/2}} |P_0(x', \psi^{-1}(y_1)) - P_0(x', \psi^{-1}(y_2))|^p \rho_0(x') \, dx' \right)$$

$$\leq c(\tilde{\mu}, \psi)(|\psi^{-1}(y_1) - \psi^{-1}(y_2)|^p + |\psi^{-1}(y_1) - \psi^{-1}(y_2)|^\gamma),$$

where $\gamma$ is as in Lemma 4.1 with $a = 1$, and we have used both assertions of this lemma to derive the last inequality. Since $\psi^{-1}$ is Lipschitz on $\overline{D}$, the bounds of the lemma follow easily. $\square$

LEMMA 4.3. *Let $p \in (1+\beta, 2)$ and $q \in (1, 1+\beta)$. There exists a constant $c = c(\beta, p, q, D, \mu)$ such that*

$$\mathbb{E}[|\theta(y_1)Z(y_1) - \theta(y_2)Z(y_2)|^q] \leq \begin{cases} c|y_1 - y_2|^q, & \text{if } 0 < \beta < \tfrac{1}{2}, \\ c|y_1 - y_2|^{q(3-p)/p}, & \text{if } \tfrac{1}{2} \leq \beta < 1. \end{cases}$$

PROOF. Recall from the proof of Theorem 1.1(ii) that the function $(s, x) \to P_D(x, y)$ belongs to $\mathbb{L}^r$ for any $y \in \partial D$ and $r \in (1, 3)$. From Lemma 3.1, it follows that

$$\mathbb{E}[|\theta(y_1)Z(y_1) - \theta(y_2)Z(y_2)|^q]$$
$$\leq c(\beta, p, q)\bigg(\bigg(\int_D |\theta(y_1)P_D(x, y_1) - \theta(y_2)P_D(x, y_2)|^p G_D\mu(x)\, dx\bigg)^{q/p}$$
$$+ \int_D |\theta(y_1)P_D(x, y_1) - \theta(y_2)P_D(x, y_2)|^q G_D\mu(x)\, dy\bigg).$$

In the case $0 < \beta < \tfrac{1}{2}$, choose $p \in (1+\beta, \tfrac{3}{2})$ and immediately get the desired bound from Lemma 4.2. Similarly, in the case $\tfrac{1}{2} \leq \beta < 1$, the desired result follows from Lemma 4.2. $\square$

PROOF OF THEOREM 1.2(a). We apply the Kolmogorov criterion of continuity to get the existence of a continuous modification of the process $\theta(y)Z(y)$ [and hence also of $Z(y)$]. The needed bounds for moments of increments of $\theta(y)Z(y)$ are obtained from the preceding lemma: In the case $0 < \beta < \tfrac{1}{2}$, this is immediate since $q > 1$, and in the case $\tfrac{1}{2} \leq \beta < 1$, we observe that we can choose $p$ and $q$ sufficiently close to $1 + \beta$ to ensure that $q(3-p)/p > 1$. The existence of a continuous modification of the process $Z(y)$, together with the remarks of the beginning of this section, completes the proof of part (a) of Theorem 1.2, in the simply connected case.

The general case when $D$ is not simply connected can be treated via a localization procedure analogous to Section 4 of [15]. Instead of the special Markov property of the Brownian snake used in [15], one uses the Markov property of superprocesses in the form stated in Theorem I.1.3 of [6]. Details are left to the reader. $\square$

**5. Irregularity of the density in high dimensions.** In this section, $3 \leq d < 1 + 2/\beta$. If $z \in \partial D$ and $r > 0$, we denote by $B_\partial(z, r)$ the open ball centered at $z$ and with radius $r$ in $\partial D$: $B_\partial(z, r) = \{y \in \partial D : |y - z| < r\}$.

In order to prove part (b) of Theorem 1.2, it is enough to verify that the property

(5.1) $\quad \|\overline{X}^D(\cdot)\|_B = \infty, \qquad \mathbb{P}$-a.s. on the event $\{X^D(B) > 0\}$

holds whenever B is a fixed boundary ball.

We thus fix a boundary ball $B = B_\partial(z_0, \eta_0)$. For technical reasons, we also introduce a smaller closed ball $B' = \overline{B}_\partial(z_0, \eta_0')$, with $\eta_0' < \eta_0$. If $\partial B'$ denotes the relative boundary of $B'$, we assume that $\sigma(\partial B') = 0$ (this is certainly true for all but countably many values of $\eta_0'$). We consider a sequence $(\varepsilon_n)$ of positive numbers decreasing to 0. For definiteness we may take $\varepsilon_n = 2^{-n}$. Then, for every integer $n \geq 1$, we set

$$B_n = \{x \in D : \operatorname{dist}(x, B') \leq \varepsilon_n\}.$$

LEMMA 5.1. *We have*

$$\varepsilon_n^{-2} \int_0^\infty X_s(B_n)\, ds \to X^D(B') \qquad \text{as } n \to \infty, \text{ in } \mathbb{P}\text{-probability.}$$

PROOF. From (1.1), we have for every $\varphi \in \mathcal{B}_b(\partial D)$,

$$\mathbb{E}[\langle X^D, \varphi \rangle] = \langle \mu, P_D \varphi \rangle = \int \mu(dx) \int_{\partial D} \sigma(dy) P_D(x, y) \varphi(y).$$

Taking $\varphi = \mathbb{1}_{\partial B'}$, we see that our assumption $\sigma(\partial B') = 0$ implies $X^D(\partial B') = 0$, a.s. The statement of the lemma is then an easy consequence of the weak convergence of $X_\varepsilon^D$ towards $X^D$ (Proposition 2.1). $\square$

We fix $\alpha \in (2/(\beta+1), 2)$. Let

$$\tau_n = \inf\{s > 0 : \Delta X_s(B_n) > \varepsilon_n^\alpha\}.$$

LEMMA 5.2. *We have*

(5.2) $\qquad \mathbb{P}(\tau_n = \infty \mid X^D(B') > 0) \to 0 \qquad \text{as } n \to \infty,$

*and*

(5.3) $\qquad \limsup_{n \to \infty} \mathbb{P}(\tau_n = \infty) \leq \mathbb{P}(X^D(B') = 0).$

PROOF. Equation (5.3) is an immediate consequence of (5.2). To verify (5.2), we will follow the lines of the proof of Lemma 4.1 of [18]. Define

$$Z_t^n = N([0,t] \times \{\mu \in \mathcal{M}_F^D : \mu(B_n) > \varepsilon_n^\alpha\}),$$

where $N$ is the point measure of jumps of the process $X$, which was introduced in Section 2.3. Then

(5.4) $\qquad \{\tau_n = \infty\} = \{Z_\infty^n = 0\}.$

Recall (2.6) for the compensator of $N$. From a classical time change result for counting processes (see, e.g., Theorem 10.33 in [13]), we get that for each $n$ there exists a standard Poisson process $A^n = (A^n(t), t \geq 0)$ such that

$$Z_t^n = A^n\Big(c(\beta)\varepsilon_n^{-\alpha(\beta+1)} \int_0^t X_s(\mathrm{B}_n)\,ds\Big),$$

where $c(\beta) = \beta/\Gamma(1-\beta) > 0$. Fix $\delta > 0$ such that $2 - \alpha(\beta+1) + \delta < 0$. Then

$$\mathbb{P}(Z_\infty^n = 0, X^D(\mathrm{B}') > 0)$$
$$\leq \mathbb{P}\Big(A^n(\varepsilon_n^{-\delta}) = 0,\ c(\beta)\varepsilon_n^{-\alpha(\beta+1)} \int_0^\infty X_s(\mathrm{B}_n)\,ds > \varepsilon_n^{-\delta},\ X^D(\mathrm{B}') > 0\Big)$$
(5.5) $\quad + \mathbb{P}\Big(c(\beta)\varepsilon_n^{-\alpha(\beta+1)} \int_0^\infty X_s(\mathrm{B}_n)\,ds \leq \varepsilon_n^{-\delta},\ X^D(\mathrm{B}') > 0\Big)$
$$\leq \mathbb{P}(A^n(\varepsilon_n^{-\delta}) = 0)$$
$$\quad + \mathbb{P}\Big(c(\beta)\varepsilon_n^{2-\alpha(\beta+1)+\delta}\Big(\varepsilon_n^{-2}\int_0^\infty X_s(\mathrm{B}_n)\,ds\Big) \leq 1,\ X^D(\mathrm{B}') > 0\Big).$$

The first term on the right-hand side of (5.5) is $\mathbb{P}(A^n(\varepsilon_n^{-\delta}) = 0) = \exp\{-\varepsilon_n^{-\delta}\}$, which converges to $0$ as $n \to \infty$. Now, by Lemma 5.1, $\varepsilon_n^{-2}\int_0^\infty X_s(\mathrm{B}_n)\,ds \to X^D(\mathrm{B}')$, in probability, as $n \to \infty$. Since $2 - \alpha(\beta+1) + \delta < 0$, we immediately get that

$$\mathbb{P}\Big(c(\beta)\varepsilon_n^{2-\alpha(\beta+1)+\delta}\Big(\varepsilon_n^{-2}\int_0^\infty X_s(\mathrm{B}_n)\,ds\Big) \leq 1,\ X^D(\mathrm{B}') > 0\Big) \to 0,$$

as $n \to \infty$. Hence, the result follows from (5.4) and (5.5). $\square$

In order to get a lower bound for $\overline{X}^D$ in terms of $X^D$, we observe that there exists a positive constant $C_1 = C_1(D)$ such that $\sigma(B_\partial(x, 2\varepsilon_n)) \leq C_1\varepsilon_n^{d-1}$ for every $n \geq 1$ and $x \in \partial D$. If $n$ is large enough so that $2\varepsilon_n < \eta_0 - \eta_0'$, which we assume from now on, we have $B_\partial(z, 2\varepsilon_n) \subset \mathrm{B}$ for every $z \in \mathrm{B}'$, and so

$$\sup_{z \in \mathrm{B}'}\langle X^D, \mathbb{1}_{B_\partial(z, 2\varepsilon_n)}\rangle \leq C_1\varepsilon_n^{d-1}\|\overline{X}^D\|_\mathrm{B}.$$

Thus

(5.6) $\quad \mathbb{E}[\exp\{-C_1\|\overline{X}^D\|_\mathrm{B}\}] \leq \mathbb{E}\Big[\exp\Big\{-\sup_{z \in \mathrm{B}'}\varepsilon_n^{1-d}\langle X^D, \mathbb{1}_{B_\partial(z, 2\varepsilon_n)}\rangle\Big\}\Big].$

On the event $\{\tau_n < \infty\}$, denote by $\zeta_n, r_n$ the spatial location and the size of the jump at time $\tau_n$, meaning that $\Delta X_{\tau_n} = r_n\delta_{\zeta_n}$. From the strong Markov property at time $\tau_n$, together with the additivity property of superprocesses, we know that conditionally on $\{\tau_n < \infty\}$, the process $(X_{\tau_n+t}, t \geq 0)$ is bounded below in distribution by $(\widetilde{X}_t^n, t \geq 0)$, where $\widetilde{X}^n$ is a super-Brownian motion with initial value $r_n \delta_{\zeta_n}$. From our approximations of the exit measure, it follows that conditionally on $\{\tau_n < \infty\}$, $X^D$ is bounded below in distribution by the exit measure $\widetilde{X}^{n,D}$ of $\widetilde{X}^n$ from $D$. Hence, from (5.6) we get

$$\mathbb{E}[\exp\{-C_1 \|\overline{X}^D\|_B\}]$$
(5.7)
$$\leq \mathbb{E}\left[\mathbb{1}_{\{\tau_n < \infty\}} \exp\left\{-\sup_{z \in B'} \varepsilon_n^{1-d} \langle X^D, \mathbb{1}_{B_\partial(z, 2\varepsilon_n)} \rangle\right\}\right] + \mathbb{P}(\tau_n = \infty)$$
$$\leq \mathbb{E}\left[\mathbb{1}_{\{\tau_n < \infty\}} \mathbb{E}_{r_n \delta_{\zeta_n}}\left[\exp\left\{-\sup_{z \in B'} \varepsilon_n^{1-d} \langle X^D, \mathbb{1}_{B_\partial(z, 2\varepsilon_n)} \rangle\right\}\right]\right] + \mathbb{P}(\tau_n = \infty).$$

Note that, on the event $\{\tau_n < \infty\}$, we have $r_n \geq \varepsilon_n^\alpha$ and $\zeta_n \in B_n$. We now claim that

(5.8) $$\lim_{n \to \infty} \sup_{x \in B_n, r \geq \varepsilon_n^\alpha} \mathbb{E}_{r \delta_x}\left[\exp\left\{-\sup_{z \in B'} \varepsilon_n^{1-d} \langle X^D, \mathbb{1}_{B_\partial(z, 2\varepsilon_n)} \rangle\right\}\right] = 0.$$

To verify (5.8), let $x_0 \in B_n$ and $r \geq \varepsilon_n^\alpha$. By the definition of $B_n$, there exists $y_0 \in B'$ such that $|y_0 - x_0| \leq \varepsilon_n$. Then, using the Laplace functional of the exit measure as recalled in Section 2.2,

(5.9)
$$\mathbb{E}_{r \delta_{x_0}}\left[\exp\left\{-\sup_{z \in B'} \varepsilon_n^{1-d} \langle X^D, \mathbb{1}_{B_\partial(z, 2\varepsilon_n)} \rangle\right\}\right]$$
$$\leq \mathbb{E}_{r \delta_{x_0}}\left[\exp\left\{-\varepsilon_n^{1-d} \langle X^D, \mathbb{1}_{B_\partial(y_0, 2\varepsilon_n)} \rangle\right\}\right]$$
$$= \exp(-r v_{y_0}^n(x_0))$$
$$\leq \exp(-\varepsilon_n^\alpha v_{y_0}^n(x_0)),$$

where the nonnegative function $(v_{y_0}^n(x), x \in D)$ solves the integral equation

(5.10) $$v_{y_0}^n(x) + \int_D G_D(x, y) v_{y_0}^n(y)^{1+\beta} \, dy = \varepsilon_n^{1-d} \int_{B_\partial(y_0, 2\varepsilon_n)} P_D(x, z) \sigma(dz).$$

LEMMA 5.3. *Under the conditions $2/(\beta + 1) < \alpha < 2$ and $3 \leq d < 1 + 2/\beta$, we have*

(5.11) $$\lim_{n \to \infty}\left(\inf_{x_0 \in B_n, y_0 \in B', |y_0 - x_0| \leq \varepsilon_n} \varepsilon_n^\alpha v_{y_0}^n(x_0)\right) = +\infty.$$

Let us postpone the proof of Lemma 5.3. Our claim (5.8) readily follows from (5.9) and (5.11). By passing to the limit $n \to \infty$ in the right-hand side of (5.7), and then using Lemma 5.2, we arrive at

$$\mathbb{E}[\exp\{-\|\overline{X}^D\|_B\}] \leq \limsup_{n \to \infty} \mathbb{P}(\tau_n = \infty) \leq \mathbb{P}(X^D(B') = 0).$$

We can now let B′ increase to B by varying $\eta_0'$ along a suitable sequence increasing to $\eta_0$. Since the event $\{X^D(\mathrm{B}) = 0\}$ is the decreasing limit of the events $\{X^D(\mathrm{B}') = 0\}$ along this sequence, we get

$$\mathbb{E}[\exp\{-\|\overline{X}^D\|_\mathrm{B}\}] \leq \mathbb{P}(X^D(\mathrm{B}) = 0).$$

Since obviously $\|\overline{X}^D\|_\mathrm{B} = 0$ on the event $\{X^D(\mathrm{B}) = 0\}$, the desired property (5.1) follows from this last bound. This completes the proof of part (b) of Theorem 1.2.

PROOF OF LEMMA 5.3. Let $n \geq 1$ and $x_0 \in \mathrm{B}_n$, $y_0 \in \mathrm{B}'$ such that $|x_0 - y_0| \leq \varepsilon_n$. In what follows we will need to assume that $n$ is sufficiently large, but our bounds will then be uniform in $x_0$ and $y_0$. To simplify notation we write $v^n = v^n_{y_0}$. Note that by (5.10), for every $x \in D$,

$$v^n(x) \leq \varepsilon_n^{1-d} \int_{B_\partial(y_0, 2\varepsilon_n)} P_D(x, z) \sigma(dz).$$

Therefore,

$$\begin{aligned}
(5.12) \quad & \int_D G_D(x_0, y) v^n(y)^{1+\beta} \, dy \\
& \leq \varepsilon_n^{(1-d)(1+\beta)} \int_D G_D(x_0, y) \left( \int_{B_\partial(y_0, 2\varepsilon_n)} P_D(y, z) \sigma(dz) \right)^{1+\beta} dy.
\end{aligned}$$

We first get a lower bound on the right-hand side of (5.10) for $x = x_0$. Since $D$ is of class $C^2$, there is a number $\alpha > 0$ such that, for every $z \in \partial D$, there exists an exterior sphere of radius $\alpha$ tangent to $\partial D$ at $z$. Suppose that $n$ is large enough so that $\varepsilon_n < \alpha$, and for $z \in \partial D$ denote by $B_n^z$ the closed ball with radius $\varepsilon_n/2$ tangent to $\partial D$ at $z$ and such that $B_n^z \cap D = \varnothing$. Then, if $x \in D$ is such that $|x - z| \leq \varepsilon_n$, the probability that a Brownian motion started at $x$ exits the domain $D$ at a point of $B_\partial(z, 2\varepsilon_n)$ is bounded below by the probability that this Brownian motion hits $B_n^z$ before exiting $B(z, 2\varepsilon_n)$. Clearly, this probability is bounded below by a constant $C_0(d) > 0$. Hence,

$$(5.13) \quad \int_{B_\partial(y_0, 2\varepsilon_n)} P_D(x_0, z) \sigma(dz) = \Pi_{x_0}(\xi_{\zeta-} \in B_\partial(y_0, 2\varepsilon_n)) \geq C_0(d).$$

We then turn to an upper bound for the integral over $D$ in the right-hand side of (5.12). It will be convenient to deal separately with the integrals over $D \cap B(y_0, \varepsilon_n^\gamma)$ and $D \cap B(y_0, \varepsilon_n^\gamma)^c$, respectively, where $0 < \gamma < 1$ is chosen so that

$$d < \frac{1+\gamma}{\beta} + 1.$$

With obvious modifications, we can then follow the calculations of ([1], page 81) and, using (2.1) and (2.3) in the first inequality below, we obtain, for $n$ large enough,

$$I_1^n \equiv \int_{D \cap B(y_0, \varepsilon_n^\gamma)^c} G_D(x_0, y) \left( \int_{B_\partial(y_0, 2\varepsilon_n)} P_D(y, z) \sigma(dz) \right)^{1+\beta} dy$$

$$\leq c(D) \left( \int_{B_\partial(y_0, 2\varepsilon_n)} \sigma(dz) \right)^{1+\beta}$$

$$\times \int_{D \cap B(y_0, \varepsilon_n^\gamma)^c} |x_0 - y|^{1-d} \rho(y)^{2+\beta} \sup_{z \in B_\partial(y_0, 2\varepsilon_n)} |y - z|^{-d(1+\beta)} dy$$

$$\leq c(D) \varepsilon_n^{(d-1)(1+\beta)} \varepsilon_n^{\gamma(1-d)} \int_{D \cap B(y_0, \varepsilon_n^\gamma)^c} \rho(y)^{2+\beta} \sup_{z \in B_\partial(y_0, 2\varepsilon_n)} |y - z|^{-d(1+\beta)} dy$$

$$\leq c(D) \varepsilon_n^{(d-1)(1+\beta-\gamma)} \int_{D \cap B(y_0, \varepsilon_n^\gamma)^c} (\mathrm{dist}(y, B_\partial(y_0, 2\varepsilon_n)))^{2+\beta-d(1+\beta)} dy$$

$$\leq c(D) \varepsilon_n^{(d-1)(1+\beta-\gamma)} \int_{\varepsilon_n^\gamma}^{\mathrm{diam}\, D} r^{d-1} (r - 2\varepsilon_n)^{2+\beta-d(1+\beta)} dr$$

$$\leq c(D) \varepsilon_n^{(d-1)(1+\beta-\gamma)},$$

where the last inequality holds because $d < 1 + 2/\beta$ implies $1 + \beta - d\beta > -1$.

Let us turn to the integral over $D \cap B(y_0, \varepsilon_n^\gamma)$, which is denoted by $I_2^n$. Notice that, for $y \in D$,

$$\int_{B_\partial(y_0, 2\varepsilon_n)} P_D(y, z) \sigma(dz) = \Pi_y(\xi_{\zeta-} \in B_\partial(y_0, 2\varepsilon_n)) \leq 1.$$

Hence, using again (2.1) and (2.3),

$$I_2^n \leq \int_{D \cap B(y_0, \varepsilon_n^\gamma)} G_D(x_0, y) \left( \int_{B_\partial(y_0, 2\varepsilon_n)} P_D(y, z) \sigma(dz) \right) dy$$

$$\leq c(D) \int_{D \cap B(y_0, \varepsilon_n^\gamma)} |x_0 - y|^{1-d} \rho(y)$$

$$\times \left( \mathbb{1}_{\{\rho(y) \leq 4\varepsilon_n\}} \right.$$

$$\left. + \mathbb{1}_{\{\rho(y) > 4\varepsilon_n\}} \int_{B_\partial(y_0, 2\varepsilon_n)} \rho(y) |y - z|^{-d} \sigma(dz) \right) dy$$

$$\leq c(D) \int_{D \cap B(y_0, \varepsilon_n^\gamma)} |x_0 - y|^{1-d} \rho(y) (\mathbb{1}_{\{\rho(y) \leq 4\varepsilon_n\}} + \mathbb{1}_{\{\rho(y) > 4\varepsilon_n\}} \rho(y)^{1-d} \varepsilon_n^{d-1}) dy$$

$$\leq c(D) \varepsilon_n \int_{D \cap B(y_0, \varepsilon_n^\gamma)} |x_0 - y|^{1-d} dy$$

$$\leq c(D) \varepsilon_n \varepsilon_n^\gamma.$$

By combining the preceding bounds, we get

$$\int_D G_D(x,y) v^n(y)^{1+\beta}\,dy \leq \varepsilon_n^{(1-d)(1+\beta)}(I_1^n + I_2^n)$$
(5.14)
$$\leq c\varepsilon_n^{(1-d)(1+\beta)}(\varepsilon_n^{(d-1)(1+\beta-\gamma)} + \varepsilon_n^{1+\gamma})$$
$$= c(\varepsilon_n^{(1-d)\gamma} + \varepsilon_n^{2+\beta+\gamma-d-d\beta}).$$

Therefore, by (5.10), (5.12)–(5.14), we have

(5.15) $$v^n(x_0) \geq C_0(d)\varepsilon_n^{1-d} - c(D)(\varepsilon_n^{(1-d)\gamma} + \varepsilon_n^{2+\beta+\gamma-d-d\beta}).$$

Hence,

(5.16) $$\varepsilon_n^\alpha v^n(x_0) \geq \varepsilon_n^{\alpha+1-d}(C_0(d) - c(D)(\varepsilon_n^{(d-1)(1-\gamma)} + \varepsilon_n^{1+\beta+\gamma-d\beta}))$$

for $n$ large enough. Since $d < \frac{1+\gamma}{\beta}+1$ and $\gamma < 1$, the expression in brackets converges to $C_0(d) > 0$ as $n \to \infty$. Moreover, since $d \geq 3$ and $\alpha < 2$, we have $\varepsilon_n^{\alpha+1-d} \to +\infty$ as $n \to \infty$, and the desired result follows. $\square$

**6. The probabilistic representation of solutions of $\frac{1}{2}\Delta u = u^{1+\beta}$.** In this section, we concentrate on the case when $D$ is the unit ball of $\mathbb{R}^d$, and we prove Theorem 1.3. Before starting the proof, let us observe that our definition of the range (which agrees with [6]) is slightly different from the one in [7] or [8]. The reason is that a superprocess is defined in [7] or [8] as the collection of its exit measures from space-time open sets. It is, however, not hard to see that both definitions give rise to the same random closed set, $\mathbb{P}_\mu$-a.s. for any $\mu \in \mathcal{M}_F^D$.

We first recall the definition of the trace of a solution following [17]. Let $u$ be a nonnegative solution of the partial differential equation

(6.1) $$\tfrac{1}{2}\Delta u = u^{1+\beta} \qquad \text{in } D.$$

We define the trace $\mathrm{tr}(u)$ of $u$ on the boundary as the pair $(K,\nu)$, where $K$ is a compact subset of $\partial D$ and $\nu$ is a Radon measure on $\partial D \setminus K$, which is determined as follows:

(i) A point $y \in \partial D$ belongs to $K$ if and only if, for every relative neighborhood $U$ of $y$ in $\partial D$,

$$\lim_{r\uparrow 1}\int_U u(rz)\sigma(dz) = \infty.$$

(ii) For every continuous function $\varphi$ on $\partial D$, with compact support contained in $\partial D \setminus K$,

$$\lim_{r\uparrow 1}\int_{\partial D} u(rz)\varphi(z)\sigma(dz) = \int_{\partial D\setminus K}\varphi(z)\nu(dz).$$

Under the condition $d < 1 + 2/\beta$, Marcus and Véron [17] proved that the mapping $u \to \operatorname{tr}(u)$ induces a one-to-one correspondence between the set of all nonnegative solutions of $\frac{1}{2}\Delta u = u^{1+\beta}$ in $D$ and the set of all pairs $(K, \nu)$, where $K$ is a compact subset of $\partial D$ and $\nu$ is a Radon measure on $\partial D \setminus K$. (In the special case $\beta = 1$, this result was obtained earlier in [15].)

Let us prove the first assertion of Theorem 1.3. If $u$ is given by (1.5), we aim at proving that $u$ solves (6.1). This is basically a consequence of the known connections between superprocesses and partial differential equations. Consider first the case when $\nu(dy) = g(y)\sigma(dy)$, where $g$ is a nonnegative continuous function on $\partial D$, with support contained in $\partial D \setminus K$. The random variable $Y$ such that $Y = +\infty$ on the event $\{\mathcal{R} \cap K \neq \varnothing\}$ and $Y = \langle X^D, g \rangle$ on $\{\mathcal{R} \cap K = \varnothing\}$ is a stochastic boundary value in the sense of [7] (see, in particular, Theorem 6.1 in [7]). Therefore the function

$$x \to -\log \mathbb{E}_{\delta_x}[\exp -Y] = -\log \mathbb{E}_{\delta_x}\left[\mathbb{1}_{\{\mathcal{R} \cap K = \varnothing\}} \exp\left(-\int g(y) X^D(dy)\right)\right]$$

solves $\frac{1}{2}\Delta u = u^{1+\beta}$ in $D$.

Coming back to the case of a general Radon measure $\nu$ on $\partial D \setminus K$, we may find a sequence of nonnegative continuous functions $g_n$, with support contained in $\partial D \setminus K$, such that

$$\lim_{n \to \infty} \int \varphi(y) g_n(y) \sigma(dy) = \int \varphi(y) \nu(dy)$$

for every $\varphi \in \mathcal{C}(D)$ with compact support contained in $\partial D \setminus K$.

LEMMA 6.1. *On the event $\{\mathcal{R} \cap K = \varnothing\}$, we have*

$$\langle X^D, g_n \rangle \to \int \overline{X}^D(y) \nu(dy)$$

*as $n \to \infty$, in $\mathbb{P}_{\delta_{x_0}}$-probability for every $x_0 \in D$.*

PROOF. Let $\varepsilon > 0$ and $K_\varepsilon = \{y \in \partial D : \operatorname{dist}(y, K) < \varepsilon\}$. Since $\mathcal{R}$ is a closed set, the event $\{\mathcal{R} \cap K = \varnothing\}$ is the union of the events $\{\mathcal{R} \cap K_\varepsilon = \varnothing\}$ over all $\varepsilon > 0$. Also, on the event $\{\mathcal{R} \cap K_\varepsilon = \varnothing\}$, it is easy to see that $X^D$ puts no mass on $K_\varepsilon$ (use Proposition 2.1) and that $\overline{X}^D(y) = 0$ a.s., for every $y \in K_\varepsilon$.

Fix $\varepsilon > 0$ and let $h_\varepsilon : \partial D \to [0,1]$ be a continuous function such that $h_\varepsilon(y) = 0$ if $y \in K_{\varepsilon/2}$ and $h_\varepsilon(y) = 1$ if $y \notin K_\varepsilon$. In view of the preceding remarks, the proof of the lemma reduces to checking that

$$\lim_{n \to \infty} \langle X^D, h_\varepsilon g_n \rangle = \int \overline{X}^D(y) h_\varepsilon(y) \nu(dy),$$

in $\mathbb{P}_{\delta_{x_0}}$-probability, for any $x_0 \in D, \varepsilon > 0$.

As a special case of (1.1), we have $\mathbb{P}_{\delta_{x_0}}$-a.s.,

$$\langle X^D, h_\varepsilon g_n \rangle = P_D(h_\varepsilon g_n)(x_0) + \int_0^\infty \int_D P_D(h_\varepsilon g_n)(x) M(ds, dx).$$

Now, for every $x \in D$, we have

$$P_D(h_\varepsilon g_n)(x) = \int_{\partial D} P_D(x,y) h_\varepsilon(y) g_n(y) \sigma(dy) \to \int_{\partial D} P_D(x,y) h_\varepsilon(y) \nu(dy),$$

as $n \to \infty$. Recall from the proof of Theorem 1.1(ii) that the (time-independent) functions $(s,x) \to P_D(x,z)$ are bounded in $\mathbb{L}^p$ when $z$ varies in $\partial D$ for any $p \in (1+\beta, (d+1)/(d-1))$. It follows that the previous convergence holds in $\mathbb{L}^p$ for any $p \in (1+\beta, (d+1)/(d-1))$. By Lemma 3.1, we conclude that $\langle X^D, h_\varepsilon g_n \rangle$ converges in $L^q(\mathbb{P}_{\delta_{x_0}})$, for every $q \in (1, 1+\beta)$, toward

$$\int_{\partial D} P_D(x_0, y) h_\varepsilon(y) \nu(dy) + \int_0^\infty \int_D \left( \int_{\partial D} P_D(x,y) h_\varepsilon(y) \nu(dy) \right) M(ds\, dx)$$
$$= \int_{\partial D} \overline{X}^D(y) h_\varepsilon(y) \nu(dy),$$

thanks to (1.2) and the "Fubini theorem" Lemma 3.2. $\square$

We come back to the proof of Theorem 1.3. For every $n \geq 1$ let

$$u_n(x) = -\log \mathbb{E}_{\delta_x}\left[ \mathbb{1}_{\{\mathcal{R} \cap K = \varnothing\}} \exp\left( -\int g_n(y) X^D(dy) \right) \right], \qquad x \in D.$$

We already saw that $u_n$ solves (6.1), and by the lemma, $u_n(x)$ converges to $u(x)$ as $n \to \infty$, for every $x \in D$. Since the set of nonnegative solutions of (6.1) is closed under pointwise convergence (see, e.g., Theorem 5.3.2 in [8]), we conclude that $u$ also solves (6.1). This completes the proof of the first part of Theorem 1.3.

In order to prove the second half of the theorem, we keep assuming that $u$ is given by (1.5) and we determine the trace of $u$. For every $n$, set $(K_n, \nu_n) = \text{tr}(u_n)$. Note that

$$u_n(x) \geq u_K(x) \equiv -\log \mathbb{P}_{\delta_x}(\mathcal{R} \cap K = \varnothing)$$

and that $u_K$ has trace $(K, 0)$. Indeed, $u_K$ is the maximal nonnegative solution of (6.1) that vanishes on $\partial D \setminus K$; see [8], Theorem 10.1.3. From the definition of the trace, it follows that $K_n \supset K$. On the other hand, set

$$u_{g_n}(x) = -\log \mathbb{E}_{\delta_x}\left[ \exp\left( -\int g_n(y) X^D(dy) \right) \right]$$

and recall that $u_{g_n}$ solves (6.1) with boundary condition $u_{|\partial D} = g_n$. From the bound

$$\left| \mathbb{E}_{\delta_x}\left[ \mathbb{1}_{\{\mathcal{R} \cap K = \varnothing\}} \exp\left( -\int g_n(y) X^D(dy) \right) \right] - \mathbb{E}_{\delta_x}\left[ \exp\left( -\int g_n(y) X^D(dy) \right) \right] \right|$$
$$\leq \mathbb{P}_{\delta_x}(\mathcal{R} \cap K \neq \varnothing)$$

and the previous observations on $u_K$, we see that $(u_n - u_{g_n})(x)$ converges to 0 as $x \to y$, for every $y \in \partial D \setminus K$. Thus $u_n$ has boundary value $g_n$ on $\partial D \setminus K$, and we conclude that $K_n = K$ and $\nu_n(dx) = g_n(x)\sigma(dx)$. Furthermore, we know from Theorem 5.6 of [17] that the convergence of $u_n$ to $u$ implies the convergence of $\operatorname{tr}(u_n)$ towards $\operatorname{tr}(u)$, in the sense of Definition 5.5 of [17], and we obtain that $\operatorname{tr}(u) = (K, \nu)$.

Finally, if $v$ is any nonnegative solution of (6.1) and $(K, \nu)$ is its trace, the solution $u$ defined by (6.1) has the same trace as $v$, and by the uniqueness theorem of [17], we must have $v = u$.

REMARK 6.1. The main contribution of [15] is a direct probabilistic proof of the special case $\beta = 1$ of Theorem 1.3. Note that the probabilistic representation of solutions in [15] looks a bit different because it is formulated in terms of excursion measures, which we did not introduce in the present work. Very probably (at least in the case $d = 2$ where the density $\overline{X}^D$ has a continuous modification), one could give a probabilistic proof of Theorem 1.3 along the lines of [15], without any reference to the results of [17]. On the other hand, this probabilistic approach remains restricted to the values $\beta \leq 1$, whereas the analytic results hold for any $\beta > 0$. For this reason, we chose to use the full strength of the results of [17] to give a short proof of the probabilistic representation (1.5). Also note that closely related results appear in the recent work of Dynkin and Kuznetsov; see, for example, Theorem 1.4 in [9].

## APPENDIX

PROOF OF LEMMA 4.1. First we will prove part (a) of the lemma. From the explicit formula (4.1) for the Poisson kernel, we have for every $x \in D_0$ and $y_1, y_2 \in \partial D_0$,

$$(A.1) \qquad |P_0(x, y_1) - P_0(x, y_2)| = \frac{1}{2\pi}(1 - |x|^2)\frac{2|x \cdot (y_1 - y_2)|}{|y_1 - x|^2|y_2 - x|^2},$$

where $u \cdot v$ stands for the usual scalar product in $\mathbb{R}^2$. Clearly, $1 - |x|^2 \leq 2\rho_0(x)$, and hence,

$$|P_0(x, y_1) - P_0(x, y_2)|^p \leq c\rho_0(x)^p \frac{|x \cdot (y_1 - y_2)|^p}{|y_1 - x|^{2p}|y_2 - x|^{2p}}.$$

Set

$$E_1 \equiv \{x \in D_0 : |y_1 - x| \vee |y_2 - x| \geq 3|y_1 - y_2|\},$$
$$E_2 \equiv \{x \in D_0 : |y_1 - x| \vee |y_2 - x| < 3|y_1 - y_2|\}.$$

If $x \in E_1$, we have plainly

(A.2) $$|y_1 - x| \wedge |y_2 - x| \geq 2|y_1 - y_2|,$$

(A.3) $$|y_1 - x| \wedge |y_2 - x| \geq \tfrac{2}{3}(|y_1 - x| \vee |y_2 - x|).$$

Also note that

(A.4) $$\begin{aligned}|x \cdot (y_1 - y_2)| &= |(x - \tfrac{1}{2}(y_1 + y_2)) \cdot (y_1 - y_2)| \\ &= \tfrac{1}{2}|(x - y_1) \cdot (y_1 - y_2) + (x - y_2) \cdot (y_1 - y_2)| \\ &\leq (|x - y_1| \vee |x - y_2|)|y_1 - y_2|.\end{aligned}$$

By combining (A.2)–(A.4), we obtain

$$\int_{E_1} \rho_0(x)^a |P_0(x, y_1) - P_0(x, y_2)|^p \, dx$$

$$\leq c \int_{E_1} \rho_0(x)^{a+p} \frac{|x \cdot (y_1 - y_2)|^p}{|y_1 - x|^{2p}|y_2 - x|^{2p}} \, dx$$

$$\leq c \left( \int_{|x-y_1| \wedge |x-y_2| > 2|y_1-y_2|} \rho_0(x)^{a+p} |y_1 - x|^{-2p} |y_2 - x|^{-p} \, dx \right) |y_1 - y_2|^p$$

$$\leq c \left( \int_{|x-y_1| > 2|y_1-y_2|} \rho_0(x)^{a+p} |y_1 - x|^{-3p} \, dx \right) |y_1 - y_2|^p$$

$$\leq c \left( \int_{2|y_1-y_2| \wedge 2}^{2} r^{1+a+p} r^{-3p} \, dr \right) |y_1 - y_2|^p$$

$$\leq \begin{cases} c|y_1 - y_2|^p, & \text{if } 0 < p < \dfrac{2+a}{2}, \\ c\left(\log_+ \dfrac{1}{|y_1 - y_2|} + 1\right)|y_1 - y_2|^p, & \text{if } p = \dfrac{2+a}{2}, \\ c|y_1 - y_2|^{2+a-p}, & \text{if } p > \dfrac{2+a}{2}. \end{cases}$$

Then consider the integral on $E_2$. If $x \in E_2$, we have by (A.4),

$$|x \cdot (y_1 - y_2)| \leq 3|y_1 - y_2|^2.$$

Also note that $|y_1 - x| \vee |y_2 - x| \geq \tfrac{1}{2}|y_1 - y_2|$. Then it follows that

$$\int_{E_2} \rho_0(x)^a |P_0(x, y_1) - P_0(x, y_2)|^p \, dx$$

$$\leq c \int_{E_2} \rho_0(x)^{a+p}(|y_1 - x|^{-2p} + |y_2 - x|^{-2p}) \, dx$$

$$\leq c \int_{|x-y_1| \vee |x-y_2| < 3|y_1-y_2|} \rho_0(x)^{a+p}(|y_1 - x| \wedge |y_2 - x|)^{-2p} \, dx$$

n/an/an/an/an/an/an/an/an/an/an/an/an/an/an/an/an/an/an/an/an/an/an/an/an/an/an/an/an/an/an/an/an/an/an/an/an/an/an/an/an/an/an/an/an/an/an/an/an/an/a

$$\leq c \int_{|x-y_1|<3|y_1-y_2|} \rho_0(x)^{a+p}|y_1-x|^{-2p}\,dx$$
$$\leq c \int_0^{3|y_1-y_2|} r^{1+a-p}\,dr$$
$$\leq c|y_1-y_2|^{2+a-p},$$

provided that $p < 2 + a$. Hence, the result of part (a) of the lemma follows by combining bounds on $E_1$ and $E_2$.

The proof of part (b) is easy. Define $b = \mathrm{dist}(B, \partial D_0)$, and recall that $b > 0$. Then from (A.1) we obtain, for every $y_1, y_2 \in \partial D_0$,

$$\sup_{x\in B}|P_0(x,y_1) - P_0(x,y_2)| \leq c \sup_{x\in B} \frac{|x||y_1-y_2|}{|y_1-x|^2|y_1-x|^2} \leq cb^{-4}|y_1-y_2|,$$

and the result follows. $\square$

**Acknowledgment.** We thank the referees for their careful reading of the first version of this work.

## REFERENCES

[1] ABRAHAM, R. and DELMAS, J. F. (2002). Some properties of the exit measure of super-Brownian motion. *Probab. Theory Related Fields* **122** 71–107. MR1883718
[2] CHUNG, K. L. and ZHAO, Z. (1995). *From Brownian Motion to Schrödinger's Equation.* Springer, Berlin. MR1329992
[3] DAWSON, D. (1993). Measure-valued Markov processes. *École d'Été de Probabilités de Saint Flour 1991. Lecture Notes in Math.* **1541** 1–260. Springer, Berlin. MR1242575
[4] DELLACHERIE, C. and MEYER, P. A. (1980). *Probabilités et Potentiel. Théorie des Martingales.* Hermann, Paris. MR566768
[5] DYNKIN, E. B. (1991). A probabilistic approach to one class of nonlinear differential equations. *Probab. Theory Related Fields* **90** 89–115. MR1109476
[6] DYNKIN, E. B. (1993). Superprocesses and partial differential equations. *Ann. Probab.* **21** 1185–1262. MR1235414
[7] DYNKIN, E. B. (1998). Stochastic boundary values and boundary singularities for solutions of the equation $Lu = u^\alpha$. *J. Funct. Anal.* **153** 147–185. MR1609265
[8] DYNKIN, E. B. (2002). *Diffusions, Superdiffusions and Partial Differential Equations.* Amer. Math. Soc., Providence, RI. MR1883198
[9] DYNKIN, E. B. and KUZNETSOV, S. E. (1998). Trace on the boundary for solutions of nonlinear differential equations. *Trans. Amer. Math. Soc.* **350** 4499–4519. MR1422602
[10] EL KAROUI, N. and ROELLY, S. (1991). Propriétés de martingales, explosion et représentation de Lévy–Khintchine d'une classe de processus de branchement à valeurs mesures. *Stochastic Process. Appl.* **38** 239–266. MR1119983
[11] FITZSIMMONS, P. J. (1992). On the martingale problem for measure-valued branching processes. In *Seminar on Stochastic Processes 1991* (E. Çinlar, K. L. Chung and M. J. Sharpe, eds.) 39–51. Birkhäuser, Boston. MR1172141

[12] GILBARG, D. and TRUDINGER, N. S. (1998). *Elliptic Partial Differential Equations of Second Order*, 2nd ed. Springer, Berlin. MR737190
[13] JACOD, J. (1979). *Calcul Stochastique et Problèmes de Martingales. Lecture Notes in Math.* **714**. Springer, Berlin. MR542115
[14] JACOD, J. and SHIRYAEV, A. N. (2003). *Limit Theorems for Stochastic Processes*, 2nd ed. Springer, Berlin. MR1943877
[15] LE GALL, J. F. (1997). A probabilistic Poisson representation for positive solutions $\Delta u = u^2$ in a planar domain. *Comm. Pure Appl. Math.* **50** 69–103. MR1423232
[16] LE GALL, J. F. (1999). *Spatial Branching Processes, Random Snakes and Partial Differential Equations*. Birkhäuser, Basel. MR1714707
[17] MARCUS, M. and VÉRON, L. (1998). The boundary trace of positive solutions of semilinear elliptic equations, I. The subcritical case. *Arch. Ration. Mech. Anal.* **144** 201–231. MR1658392
[18] MYTNIK, L. and PERKINS, E. (2003). Regularity and irregularity of $\beta$-stable super-Brownian motion. *Ann. Probab.* **31** 1413–1440. MR1989438
[19] PERKINS, E. (2002). Dawson–Watanabe superprocesses and measure-valued diffusions. *Ecole d'Eté de Probabilités de Saint-Flour 1999. Lecture Notes in Math.* **1781** 125–329. Springer, Berlin. MR1915445
[20] POMMERENKE, C. (1992). *Boundary Behavior of Conformal Maps*. Springer, Berlin. MR1217706
[21] SHEU, Y. (1996). On states of exit measure for superdiffusions. *Ann. Probab.* **24** 268–279. MR1387635
[22] WIDMAN, K. O. (1967). Inequalities for the Green function and boundary continuity of the gradient of solutions of elliptic differential equations. *Math. Scand.* **21** 13–67. MR239264

DÉPARTEMENT MATHÉMATIQUES
ET APPLICATIONS
ECOLE NORMALE SUPÉRIEURE
45 RUE D'ULM
75005 PARIS
FRANCE
E-MAIL: legall@dma.ens.fr

FACULTY OF INDUSTRIAL ENGINEERING
AND MANAGEMENT
TECHNION—ISRAEL INSTITUTE OF TECHNOLOGY
HAIFA 32000
ISRAEL
E-MAIL: leonid@ie.technion.ac.il